\documentclass[fleqn]{mat01}
\usepackage{times,mathtimy,amssymb}
\begin{document}

\setcounter{page}{309}
\firstpage{309}

\newcommand{\rr}{\mathbb{R}^d}
\newcommand{\zz}{\mathbb{Z}}
\newcommand{\pp}{\mathbb{P}}
\newcommand{\ee}{\mathbb{E}}
\newcommand{\nn}{\mathbb{N}}
\renewcommand{\div}{\mathrm{div}}
\renewcommand{\phi}{\varphi}
\renewcommand{\epsilon}{\varepsilon}
\renewcommand{\kappa}{\varkappa}
\newcommand{\<}{\langle}
\renewcommand{\>}{\rangle}
\renewcommand{\le}{\leqslant}
\renewcommand{\ge}{\geqslant}

\newtheorem{theore}{Theorem}
\renewcommand\thetheore{\arabic{section}.\arabic{theore}}
\newtheorem{theor}[theore]{\bf Theorem}
\newtheorem{rem}[theore]{Remark}

\newtheorem{propo}[theore]{\rm PROPOSITION}
\newtheorem{lem}[theore]{Lemma}
\newtheorem{definit}[theore]{\rm DEFINITION}
\newtheorem{coro}[theore]{\rm COROLLARY}
\newtheorem{exampl}[theore]{Example}

\title{Random walks in a random environment\footnote{This is essentially a transcript of the plenary talk given at the Joint India--AMS Mathematics Meeting held in December 2003 in Bangalore, India.}}

\markboth{S~R~S~Varadhan}{Random walks in a random environment}

\author{S~R~S~VARADHAN}

\address{Department of Mathematics, Courant Institute of Mathematical
Sciences, New York University, NY 10012, USA\\
\noindent E-mail: varadhan@cims.nyu.edu}

\volume{114}

\mon{November}

\parts{4}

\Date{MS received 16 June 2004}

\begin{abstract}
Random walks as well as diffusions in random media are considered.
Methods are developed that allow one to establish large deviation
results for both the `quenched' and the `averaged' case.
\end{abstract}


\maketitle

\section{Introduction}

A random walk on ${\zz}^d$ is a stochastic process $\{S_n\hbox{:}\ n\ge
0\}$ defined as 
\begin{equation*}
S_0=0,\quad S_n=X_1+\cdots+X_n \quad {\rm  for}\  n\ge 1,
\end{equation*} 
where $\{X_i\}$ are independent and identically distributed random
variables with a common distribution $\pi(z)=P[X_i=z]$ for $z\in
{\zz}^d$. The relationship between the behavior of $S_n$ as $n\to\infty$
and properties of $\pi(\cdot)$ is quite well-understood and is easy to
analyse, using often Fourier analysis as the main tool.

For instance if $\sum_{z\in{\zz}^d} |z|\,\pi(z)<\infty$, the law of
large numbers states that
\begin{equation*}
\lim_{n\to\infty}\frac{S_n}{n}=\sum_z z\,\pi(z).
\end{equation*}

If $\sum_{z\in{\zz}^d} |z|^2\,\pi(z)<\infty,$ and $\sum_{z\in{\zz}^d
}z\,\pi(z)=0$, then the central limit theorem asserts that the
distribution of ${S_n}/{\sqrt n}$ is asymptotically Gaussian on
$\rr$ with mean $0$ and covariance $C=\{C_{r,s}\}$ given by
\begin{equation*}
C_{r,s}=\sum_{z\in{\zz}^d} z_rz_s\,\pi(z).
\end{equation*}

Cram\'{e}r's theorem on large deviations provides exponential rates of
convergence of 
\begin{equation*}
P \left[\frac{S_n}{n}\in A \right]
\end{equation*}
for sets $A$ excluding the mean $m=\sum_{z\in{\zz}^d} z\, \pi(z)$.
Assuming that 
\begin{equation*}
M(\theta)=\sum_{z\in{\zz}^d} {\rm e}^{\<z,\theta\> }\pi(z)<\infty
\end{equation*}
for $\theta\in\rr$, for a wide class of sets $A$,
\begin{equation*}
\lim_{n\to\infty}\frac{1}{n}\log P \left[\frac{S_n}{n}\in A \right]
\end{equation*}
exists and is given by $-\inf_{x\in A} I(x)$ where
\begin{equation*}
I(x)=\sup_{\theta\in \rr}[\<\theta, x\>-\log M(\theta)].
\end{equation*}

The situation is quite different if we drop the hypothesis of
independence of increments of $S_n$ and replace it with the assumption
that $S_n$ is a Markov process. In other words, there exists $p( z', z)$
such that
\begin{equation*}
P[S_{n+1}=z| S_1, S_2,\ldots, S_n]=p(S_n, z).
\end{equation*}
The random walk case considered earlier is a special case where $p(z',
z)$ takes the form $\pi(z-z')$. Sometimes it is more convenient to
denote $p(z',z)$ by $\pi (z', z-z')$ so that the independent case is
when $\pi(z', \cdot)$ is independent of $z'$. In this generality not
much can be said. However if $\pi(z',z)$ is periodic in $z'$, the phase
gets averaged out and some homogenization takes place with $\pi (z', z)$
getting effectively replaced by an averaged $\hat\pi(\cdot)$. Results on
law of large numbers, central limit theorem and large deviation
estimates can be established, with some extra work, quite similar to the
independent case.

We will investigate the case when $\pi(z',z)$ are random but spatially
homogeneous in $z'$. This can viewed as a limiting version of the
periodic case when the period becomes large.

\section{Formulation}

We will start with a probability space $(\Omega, \Sigma, P)$ on which
${\zz}^d$ acts ergodically as a family $\tau_z$ of measure preserving
transformations. We are given $\pi(\omega, z)$ which is a probability
distribution on ${\zz}^d$ for each $\omega$, that are measurable
functions of $\omega$. One can then generate random transition
probabilities $\pi(\omega, z',z)$ by defining
\begin{equation*}
\pi(\omega, z',z)=\pi(\tau_{z'}\omega,z).
\end{equation*}
For each $\omega$, $\pi(\omega, z',z)$ can serve as the transition
probability of a Markov process on ${\zz}^d$ and the measure
corresponding to this process, starting from $0$, is denoted by
$Q^\omega$. This is of course random (depends on $\omega$) and is called
the random walk in the random environment $\omega$. One can ask the same
questions about this random walk and in some form they may be true for
almost all omega with respect to $P$. The law of large numbers, if
valid, will take the form
\begin{equation}\label{lln}
P \left[ \,\omega\hbox{:}\ \lim_{n\to\infty}\frac{S_n}{n}=m(P) \quad
{\rm a.e.}\ Q^\omega \right] = 1.
\end{equation}
Such statement concerning the almost sure behavior of $Q^\omega$ for
almost all $\omega$ with respect to $P$ are said to deal with the \lq
quenched\rq\, version. Sometimes one wishes to study the behavior of
the \lq averaged\rq \, measure
\begin{equation}\label{average}
\bar{Q} = \int Q^\omega P({\rm d}\omega).
\end{equation}
The law of large numbers is the same because (\ref{lln}) is equivalent
to
\begin{equation*}
\bar{Q} \left[\,\omega\hbox{:}\ \lim_{n\to\infty}\frac{S_n}{n}=m(P)
\right] = 1.
\end{equation*}

On the other hand, questions on the asymptotic behavior of probabilities,
like for instance, the central limit theorem or large deviations could
be different for the quenched and the averaged cases.

A special environment called the \lq product environment\rq\, is one in
which $\pi(\omega, z',z)$ are independent for different $z'$ and have a
common distribution $\beta$ which is a probability measure on the space
$\mathcal M$ of all probability measures on ${\zz}^d$. In this case the
canonical choice for $\Omega$ is the countable product of $\mathcal M$.
The product measure having marginals $\beta$ is of course the choice for
$P$.

One can consider the continuous versions of these problems. For instance
instead of the action of ${\zz}^d$ we can have $\rr$ acting on $(\omega,
\Sigma, P)$ ergodically and consider a diffusion on $\rr$ with a random
infinitesimal generator
\begin{equation}\label{randomdiff}
({\mathcal L}^\omega u)(x) = \frac{1}{2}(\Delta u)(x) + \<b(\omega, x),
(\nabla u)(x)\>
\end{equation}
acting on smooth function on $\rr$. Here $b(\omega,x)$ is generated from
a map $b(\omega)\hbox{:}\ \Omega\to\rr$ by the action $\{\tau_x\}$ of $\rr$,
\begin{equation*}
b(\omega,x)=b(\tau_x\omega).
\end{equation*}
Again there is the quenched measure $Q^\omega$ that corresponds to the diffusion with generator ${\mathcal L}^\omega$ that starts from $0$ at time $0$ and the averaged measure that is given by the same formula (\ref{average}). This model is referred to as diffusion with a random drift. Exactly the same questions can be asked in this context.

\section{Results: One-dimensional case}

We shall describe the results in the one-dimensional case, which is
quite well-understood. A very good reference with a full bibliography is
a survey by Zeitouni \cite{Z1} at the Beijing congress as well as his
notes from St. Flour \cite{Z2}. For simplicity we will take the special
case of the product environment in which only $\pi (\omega, z, z\pm 1)$
are non-zero and are random variables $p(z)$, $q(z)=1-q(z)$ that are
independent for different values of $z$, $\beta$ being the distribution
on $[0,1]$ of $p(0)$. To avoid some technical details we shall assume
that $\beta$ has no mass near the edges i.e. is supported on a closed
subinterval of $(0,1)$.

The question of when $S_n$ is transient was answered in \cite{S}.

\begin{theor}[\!]
In order that
\begin{equation*}
Q^\omega \left[\lim_{n\to\infty} S_n=+\infty\,\right]=1
\end{equation*}
a.e. $P$ it is necessary and sufficient that
\begin{equation}\label{rec.cond}
E \left[\log\frac{p(0)}{q(0)} \right] = \int\log\frac{p}{1-p} \beta({\rm
d}p) >0.
\end{equation}
\end{theor}

In the original random walk where $p(z)\equiv p$ is a constant this is
equivalent to $p>\frac{1}{2}$ or $p-q>0$. In such a case the law of
large numbers asserts that
\begin{equation*}
\lim_{n\to\infty}\frac{S_n}{n}=p-q>0.
\end{equation*}
In our context the law of large numbers, also proved in \cite{S}, states
the following theorem.

\begin{theor}[\!]\label{lln2}
\begin{equation*}
\lim_{n\to\infty}\frac{S_n}{n}=m
\end{equation*}
a.e. $Q^\omega$ for almost all $\omega$ with respect to $P$ where $m$ is
given by
\begin{equation}
m = \begin{cases}
\frac {1-E[\frac{q}{p}]}{1+E[\frac{q}{P}]}, &if\  E[\frac{q}{p}]<1,\\[.3pc]
0, &otherwise.
\end{cases}
\end{equation}
\end{theor}

It is therefore possible that $E[\log\frac{p}{q}]>0$ while $
E[\frac{q}{p}]\ge 1$, in which case $S_n\to\infty$ but $\frac{S_n}{n}\to
0$. In the original random walk this could not happen. In the 
one-dimensional case further detailed analysis can be carried out in this
case. The problem comes from \lq traps\rq\, where the environment
conspires to hold the particle for a long time. As the particle moves to
$\infty$, it encounters increasingly deeper traps that hold it even
longer producing a critical slowing down. This phenomenon discovered by
Sinai \cite{Sinai} is well-understood. In one-dimension the traps can not be avoided.
Therefore the one-dimensional case is strikingly different from the
higher dimensional ones. 

Assuming that $m>0$, one can study the behavior of
$\xi_n = {(S_n-nm)}/{\sqrt n}$ and ask for a central limit theorem. It
turns out that the fluctuations of $\xi_n$ have two sources; from the
environment $\{p(x)\hbox{:}\ 0\le x\le nm\}$ and fluctuations of the walk in
this environment. Under $Q^\omega$, the distribution of $\xi_n$ is
Gaussian with a random mean which depends on the environment and has an
asymptotic Gaussian distribution by itself. Hence while there is a
central limit theorem for the averaged $\bar{Q}$ there is none for
$Q^\omega$. This phenomenon is not expected to persist in higher
dimensions.

There are large deviation results regarding the limits
\begin{equation*}
\lim_{n\to\infty}\frac{1}{n}\log Q^\omega \left[\frac{S_n}{n}\simeq
a\right] = I(a)
\end{equation*}
and
\begin{equation*}
\lim_{n\to\infty} \frac{1}{n}\log \bar{Q} \left[\frac{S_n}{n}\simeq
a \right] = \bar{I} (a).
\end{equation*}

The difference between $I$ and $\bar{I}$ has a natural explanation
in terms of the large deviation behavior of the environment. It is
related to the following question. If at time $n$, a particle is near
$na$ instead of being around $nm$ which is what we should expect,
what does it say about the environment in $[0, na]$? Did the particle
behave strangely in a normal environment or did it encounter a strange
environment? It is in fact a combination of both as is seen in the
relation of $\bar{I}$ to $I$ which is essentially Bayes' rule. See
for instance \cite{CGZ} and \cite{GH}.

\section{Higher dimension}

If $d\ge 2$ the law of large numbers and central limit theorems are not
that well-understood. For instance the questions of recurrence,
transience and the existence of a non-zero limit
$m=\lim_{n\to\infty}\frac{S_n}{n}$ have only partial answers with
sufficient conditions that are not that easy to check. Even when we
limit ourselves to a product environment there is no clear analog of the
condition (\ref{rec.cond}). 

However if $d\ge 3$, under some symmetry conditions that ensure that
$m=0$, a central limit theorem was proved by Bricmont and Kupiainen
\cite{BK}, provided the randomness in the environment is small. A
recent preprint by Sznitman and Zeitouni \cite{SZ} provides an alternate
proof.

Large deviation results however exist in much wider generality, both in
the quenched and the averaged cases. The large deviation principle is
essentially the existence of the limits
\begin{equation*}
\lim_{n\to\infty}\frac{1}{n}\log E[\exp[\<\theta, S_n\>]=\Psi(\theta).
\end{equation*}
The expectation is with respect to $Q^\omega$ or $\bar{Q}$, which
could produce different limits for $\Psi$. The law of large numbers and
central limit theorem involve the differentiability of $\Psi$ at
$\theta=0$, which is harder. In fact large deviation results have been
proved by general sub-additivity arguments for the quenched case.
Roughly speaking fixing $\omega$,
\begin{equation*}
Q^\omega[S_{k+\ell}\simeq (k+\ell) a]\ge Q^\omega[S_k\simeq ka]\times
Q^{\tau_{ka}\omega}[S_\ell\simeq \ell a].
\end{equation*}
One can then use some version of the sub-additive ergodic theorem. This
idea has to be cleaned up a bit, but is not hard. See for instance
\cite{Zerner} or \cite{V}. We will however explore, in the next two
sections a proof of the large deviation principle for the averaged case
and an alternate proof for the quenched case.

\section{Large deviations: The averaged case}

We want to prove that the limits
\begin{equation}\label{ldpa}
\lim_{n\to\infty} \frac{1}{n}\log \bar{Q} \left[\frac{S_n}{n}\simeq
a\right] = -\bar{I} (a)
\end{equation}
or equivalently
\begin{equation*}
\lim_{n\to\infty} \frac{1}{n} \log E^{ \bar{Q}} [\exp[\<\theta,
S_n\>]] = \bar{\Psi} (\theta)
\end{equation*}
exist. The problem is that the measure $\bar{Q}$ is not very nice.
As the random walk explores ${\zz}^d$ it learns about the environment
and in the case of the product environment, when it returns to a site
that it has visited before, the experience has not been forgotten and
leads to long term correlations. However, if we are interested in the
behavior $S_n\simeq na$ with $a\not=0$, the same site is not visited too
often and the correlations should decay fast. This can be exploited to
provide a proof of (\ref{ldpa}).

One can use Bayes' rule to calculate the conditional distribution 
\begin{equation*}
\bar{Q} [S_{n+1}=S_n+z| S_1, S_2,\ldots, S_n]=q(z|w),
\end{equation*}
where $w$ is the past history of the walk. Before we do that it is more
convenient to shift the origin as we go along so that the current
position of the random walk is always the origin and the current time is
always $0$. Then an $n$ step walk looks like $w=\{S_0=0, S_{-1}, \ldots,
S_{-n}\}$. We pick a $z$ with probability $q(z|w)$. We get a new walk of
$n+1$ steps $w'=\{S'_0=0, S'_{-1},\ldots, S'_{-(n+1)}\}$ given by
$S'_{-(k+1)}=S_{-k}-z$ for $k\ge 0$. We can now calculate $q(z|w)$. We
need to know all the numbers $\{k(w,x,z)\}$ of the number of times the
walk has visited $x$ in the past and jumped from $x$ to $x+z$. It is not
hard to see that the {\it a posteriori} probability can be calculated
as
\begin{equation*}
q(z|w)=\frac{\int \pi(z) \Pi_{z'} \pi(z')^{k(w, 0,z')} \beta({\rm
d}\pi)}{\int \Pi_{z'} \pi(z')^{k(w, 0,z')} \beta({\rm d}\pi)}.
\end{equation*}
While this makes sense initially only for walks of finite length it can
clearly be extended to all transient paths. Note that although we only
use $k(w,0,z')$, in order to obtain the new $k(w', 0,z')$ we would need
to know $k(w,z,z')$.

Let us suppose that $R$ is a process with stationary increments
$\{z_j\}$. If the increments process were ergodic and had a non-zero mean
$a$ we can again make the current position the origin and the process
will be transient. $ q(z|w)$ would exist a.e. $R$ and can be compared to
the corresponding conditional probabilities $r(z|w)$ under $R$. The
relative entropy
\begin{equation*}
H(R)=E^R \left[\sum_z r(z|w)\log\frac{r(z|w)}{q(z|w)}\right]
\end{equation*}
is then well-defined.

\setcounter{theore}{0}
\begin{theor}[\!]
The function
\begin{equation*}
\bar{I} (a) = \inf\limits_{\substack{R{\rm :}\ \int z_1\ {\rm
d}R=a\\ R\ {\rm ergodic}}} H(R)
\end{equation*}
defined for $a\not=0$ extends as a convex function to all of $R^d${\rm ,}
and with this $\bar{I}${\rm , (\ref{ldpa})} is satisfied.
\end{theor}

The proof, which uses basic large deviation techniques, can be found in
\cite{V}.

\setcounter{theore}{0}
\begin{rem}{\rm
It is interesting to note that the rate functions $I(0)$ and 
$\bar{I}(0)$ coincide at $0$ and can be calculated explicitly. Let $\cal C$ be
the convex hull of the support of $\beta$ in ${\mathcal M}({\zz}^d)$.
Then
\begin{equation*}
I(0) = \bar{I}(0) = -\inf_{\pi\in\mathcal C} \inf_{\theta\in \rr}
\log\sum_z {\rm e}^{\<\theta, z\>} \pi(z).
\end{equation*}}
\end{rem}

In particular $I(0)>0$ if and only if $0$ is not in the range of $\sum
z\pi(z)$ as $\pi$ varies over $\mathcal C$. This is referred to in the
literature as the \lq non-nestling\rq\, case.

\section{Large deviations: The quenched case}

Although a proof using the subadditive ergodic theorem exists, we will
provide an alternate approach that is more appealing. We will illustrate
this in the context of Brownian motion with a random drift
(\ref{randomdiff}).

We can define a diffusion on $\Omega$ with  generator
\begin{equation*}
{\mathcal L}=\frac{1}{2}\Delta +\<b(\omega),\nabla\>,
\end{equation*}
where $\nabla=\{D_i\}$ are the generators of the translation group
$\{\tau_x\hbox{:}\ x\in \rr\}$. This is essentially the image of lifting the
paths $x(t)$ of the diffusion on $\rr$ corresponding to ${\mathcal
L}^\omega$ to $\Omega$ by
\begin{equation*}
\omega(t)=\tau_{x(t)}\omega.
\end{equation*}
While there is no possibility of having an invariant measure on $\rr$,
on $\Omega$ one can hope to find an invariant density $\phi(\omega)$
provided we can find $\phi(\omega)\ge 0$ in $L_1(P)$, that solves
\begin{equation*}
\frac{1}{2}\Delta\phi=\nabla\cdot (b\phi).
\end{equation*}
If such a $\phi$ exists, then we have an ergodic theorem for the
diffusion process $Q^\omega$ corresponding $\mathcal L$ on $\Omega$,
\begin{equation}\label{ergodicth}
\lim_{t\to\infty}\frac{1}{t}\int_0^t f(\omega(s)) {\rm d}s = \int
f(\omega)\phi(\omega) {\rm d}P\quad {\rm a.e.}\ Q^\omega\quad {\rm
a.e.}\ P.
\end{equation}
This translates to an ergodic theorem on $R^d$  as well
\begin{equation*}
\lim_{t\to\infty}\frac{1}{t}\int_0^t f(\omega, x(s)) {\rm d}s = \int
f(\omega)\phi(\omega) {\rm d}P\quad {\rm a.e.}\ Q^\omega\quad {\rm a.e.}\ P
\end{equation*}
where now $Q^\omega$ is the quenched process in the random environment.
Since
\begin{equation*}
x(t) = \int_0^t b(\omega, x(s) {\rm d}s + \beta(t),
\end{equation*}
it is clear that
\begin{equation*}
\lim_{t\to\infty}\frac{x(t)}{t}=\int b(\omega)\phi(\omega) {\rm
d}P\quad {\rm a.e.}\ Q^\omega\quad {\rm a.e.}\ P
\end{equation*}
providing a law of large numbers for $x(t)$. While we can not be sure of
finding $\phi$ for a given $b$ it is easy to find a $b$ for a given $\phi$.
For instance, we could take $b= {\nabla\phi}/{2\phi}$. Or more
generally $b= ({\nabla\phi}/{2\phi}) + ({c}/{\phi})$ with $\nabla\cdot
c=0$. If we change $b$ to $b'= ({\nabla\phi}/{2\phi}) + c$ with
$\nabla\cdot c=0$, the new process will have relative entropy
\begin{equation*}
E^{Q^{b',\omega}} \left[\frac{1}{2}\int_0^t
\left\|b(\omega(s))-\frac{\nabla\phi(\omega(s))}{2\phi(\omega(s))}-
\frac{c(\omega(s))}{\phi(\omega(s))} \right\|^2 {\rm d}s\right].
\end{equation*}

Moreover, for almost all $\omega$ with respect to $P$, almost surely
with respect to $Q^{b',\omega}$,
\begin{equation*}
\lim_{t\to\infty} \frac{x(t)}{t} = \int \left[ \frac{\nabla\phi}{2\phi}
+ \frac{c}{\phi} \right]\phi {\rm d}p=\int c {\rm d}P.
\end{equation*}
If we fix $\int c {\rm d}P = a$, the bound
\begin{equation*}
\liminf_{t\to\infty} \frac{1}{t}\log Q^\omega \left[ \frac{x(t)}{t}
\simeq a \right] \ge -\frac{1}{2} \int \left\| b -
\frac{\nabla\phi}{2\phi} - \frac{c}{\phi} \right\|^2 \phi {\rm d}P
\end{equation*}
is easily obtained. If we define 
\begin{equation*}
I(a) = \inf\limits_{\substack{\nabla\cdot c=0\\ \int c {\rm d}P=a}}
\frac{1}{2} \int \left\| b - \frac{\nabla\phi}{2\phi} - \frac{c}{\phi}
\right\|^2 \phi {\rm d}P,
\end{equation*}
then
\begin{equation*}
\liminf_{t\to\infty} \frac{1}{t} \log Q^\omega \left[ \frac{x(t)}{t}
\simeq a \right] \ge - I(a).
\end{equation*}
Of course these statements are valid a.e. $Q^\omega$ a.e. $P$. One can
check that $I$ is convex and the upper bound amounts to proving the dual
estimate
\begin{equation*}
\lim_{t\to\infty} \frac{1}{t} \log E^{Q^\omega} [{\rm e}^{\<\theta,
x(t)\>} ]\le \Psi(\theta),
\end{equation*}
where
\begin{equation*}
\psi(\theta) = \sup_a [\<a,\theta\>-I(a)].
\end{equation*}
We need a bound on the solution of
\begin{equation*}
u_t = \frac{1}{2}\Delta u + \<b,\nabla u\>
\end{equation*}
with $u(0) = \exp[\<\theta, x\>]$. By Hopf-Cole transformation $v=\log
u$ this reduces to estimating
\begin{equation*}
v_t = \frac{1}{2}\Delta v +\frac{1}{2}\|\nabla v\|^2+\< b,\nabla v\>
\end{equation*}
with $v(0)=\<\theta, x\>$. This can be done if we can construct a
subsolution
\begin{equation*}
\frac{1}{2}\nabla\cdot w+\frac{1}{2}\|\nabla w\|^2+\< b, w\>\le
\psi(\theta)
\end{equation*}
on $\Omega$, where $w\hbox{:}\ \Omega\to\rr$ satisfies $\int w {\rm
d}P=\theta$ and $w$ is closed in the sense that $D_iw_j=D_jw_i$. The
existence of the subsolution comes from convex analysis.
\begin{align*}
\psi(\theta) &= \sup_{\substack{\phi\\ \nabla\cdot c=0}}\left[\int \<c,\theta\> {\rm d}P -\frac{1}{2}\int \left\|b-\frac{\nabla\phi}{2\phi}-\frac{c}{\phi} \right\|^2 \phi {\rm d}P\right]\\[.3pc]
&= \sup_{\phi}\sup_c \inf_{u}\left[\int \<c,\theta+\nabla u \> {\rm d}P -\frac{1}{2}\int \left\|b-\frac{\nabla\phi}{2\phi}-\frac{c}{\phi} \right\|^2 \phi {\rm d}P \right]\\[.3pc]
&= \sup_{\phi} \inf_{u}\sup_c\left[\int \<c,\theta+\nabla u\> {\rm d}P -\frac{1}{2}\int \left\|b-\frac{\nabla\phi}{2\phi}-\frac{c}{\phi}\right\|^2 \phi {\rm d}P\right]\\[.3pc]
&= \sup_{\phi} \inf_{u}\int  \sup_c\left[\<c,\theta+\nabla u \> -\frac{1}{2} \left\|b-\frac{\nabla\phi}{2\phi}-\frac{c}{\phi}\right\|^2 \phi\right] {\rm d}P%\\[.3pc]
\end{align*}
\begin{align*}
&= \sup_{\phi} \inf_{u}\int  \left[ \left\<b-\frac{\nabla \phi}{2\phi},\theta +\nabla u \right\>+\frac{1}{2}\left\|\theta+\nabla u\right\|^2\right]\phi {\rm d}P\\[.3pc]
&= \sup_{\phi} \inf_{u}\int \left[ \left[\<b,\theta +\nabla u\>+\frac{1}{2}\left\|\theta+\nabla u\right\|^2\right]\phi -\frac{1}{2}\<\nabla u,\nabla \phi\>\right]{\rm d}P\\[.3pc]
&= \sup_{\phi} \inf_{u}\int  \left[\frac{1}{2}\Delta u +\<b,\theta +\nabla u\>+\frac{1}{2}\left\|\theta+\nabla u\right\|^2\right]\phi {\rm d}P\\[.3pc]
&= \sup_{\phi} \inf_{\substack{w\ {\rm  closed}\\ \int w {\rm d}P=\theta}}\int  \left[\frac{1}{2}\nabla\cdot w+ \<b,w\>+\frac{1}{2}\|w\|^2\right]\phi {\rm d}P\\[.3pc]
&= \inf_{\substack{w\ {\rm  closed}\\ \int w {\rm d}P=\theta}}\sup_{\phi}\int  \left[\frac{1}{2}\nabla\cdot w+ \<b,w\>+\frac{1}{2}\|w\|^2\right]\phi {\rm d}P\\[.3pc]
&= \inf_{\substack{w\ {\rm  closed}\\ \int w {\rm d}P=\theta}}\sup_{\omega} \left[\frac{1}{2}\nabla\cdot w+ \<b,w\>+\frac{1}{2}\|w\|^2\right]
\end{align*}
which proves the existence of the subsolution.

\setcounter{theore}{0}
\begin{rem}{\rm 
This can be viewed as showing the existence of a limit as $\epsilon\to
0$ (homogenization) of the solution of
\begin{equation*}
u^\epsilon_t = \frac{\epsilon}{2}\Delta u^\epsilon +\frac{1}{2}\|\nabla
u^\epsilon\|^2+ \left\< b \left( \frac{x}{\epsilon}, \omega \right),
\nabla u^\epsilon \right\>
\end{equation*}
with $u^\epsilon(0, x)=f(x)$. The limit satisfies
\begin{equation*}
u_t=\Psi(\nabla u)
\end{equation*}
with $u(0,x)=f(x)$.}
\end{rem}

This can be generalized to equations of the form
\begin{equation*}
u^\epsilon_t=\frac{\epsilon}{2}\Delta u^\epsilon + H\left(
\frac{x}{\epsilon}, \nabla u^\epsilon, \omega \right).
\end{equation*}

We have left out the details as well as the hypothesis needed to prove
the results. They will appear in \cite{KRV}.

\section*{Acknowledgments}

This research was supported by a grant from the National Science
Foundation DMS-0104343.

\end{document}